# Four Fundamental Questions in Probability Theory and Statistics

**Paolo Rocchi**
IBM, via Shangai 53, Roma, Italy
LUISS University, via Alberoni 7, Roma, Italy

**Abstract** – This study has the purpose of addressing four questions that lie at the base of the probability theory and statistics, and includes two main steps. As first, we conduct the textual analysis of the most significant works written by eminent probability theorists. The textual analysis turns out to be a rather innovative method of study in this domain, and shows how the sampled writers – no matter he is a frequentist or a subjectivist – share a similar approach. Each author argues on the multifold aspects of probability then he establishes the mathematical theory on the basis of his intellectual conclusions. It may be said that mathematics ranks second.
Hilbert foresees an approach far different from that used by the sampled authors. He proposes to axiomatize the probability calculus notably to describe the probability concepts using purely mathematical criteria. In the second stage of the present research we address the four issues of the probability theory and statistics following Hilberts' recommendations. Specifically, we use two theorems that prove how the frequentist and the subjectivist models are not incompatible as many believe. Probability has distinct meanings under different hypotheses, and in turn classical statistics and Bayesian statistics are available for adoption in different circumstances. Subsequently, these results are commented upon, followed by our conclusions

**Keyword**s: Foundational issues; frequentist and subjective probabilities; classical and Bayesian statistics.

## 1. INTRODUCTION

Pascal inaugurated the modern calculus of probability in 1654. Bernoulli, De Moivre, Leibniz and others addressed complex problems by providing complete solutions, but the first attempt to define the concept of probability was authored by Laplace in 1812, over one century and half after Pascal. The Laplacian definition had non-trivial weak points, statistics began to gain weight in science and economics during the 19th century and many people, whether erudite or illiterate, felt the need for the precise description of the probability concept. Significant attempts have been made to answer this issue but theorists progressed slowly [1]. We make the historical list of the constructions that the current literature recognizes as the most telling [2]:

- The *frequentist* theory by von Mises in 1928,
- The *subjective* theory by Ramsey and de Finetti (independently) in 1931,
- The *axiomatic* theory by Kolmogorov in 1933,
- The *Bayesian* theory by Savage in 1954,
- The *logical* theory by Carnap in 1962.

We mean to focus on the most popular models of probability, the *frequentist* and the *subjective*, as long as authors agree in merging the subjective and Bayesian views into a single box. Frequentists define probability as the limit of the relative frequency in a large number of trials. Subjectivists see probability as an individual person's measure of belief that an event will occur. These interpretations show apparent incongruities and the concept of probability still emerges as a conundrum three hundred and fifty years after Pascal.

## 2. A TEXTUAL ANALYSIS

The dispute raised between the frequentist and subjective schools attracted the attention of several theorists. Papers and books have been filled with annotations and sharp comments about this complex argument. We share the same interest but decided to follow a different direction. Instead of formulating remarks and explanations, we have conducted a textual analysis of the principal *book* (B) and *essay* (E) prepared by each one of the following seven writers: John Venn, Richard von Mises, Hans Reichenbach, John Maynard Keynes, Frank Plumpton Ramsey, Bruno de Finetti, and Leonard Jimmie Savage (Table 1). We have selected these authors because they are recognized as the founders of the frequentist and subjective schools. We have overlooked other researchers in the fields who propose interesting solutions but in substance belong to one of the two schools. For example, Cox defines probability as a measure of a degree of belief that is consistent with Boolean logic, and Knuth further generalizes this scheme to include other algebras and hence other problems in science and mathematics. However, Cox and Knuth take a tacit step in advance since they adhere to the subjectivist circle.

| Work | | | | Criticism | | | |
|---|---|---|---|---|---|---|---|
| | Author | Reference | Type | Extent W (Pages) | Chapter Id | Location | Extent C (Pages) | C/W |
| | Kolmogorov | [3] | B | 84 | None | N.A. | 0 | N.A. |
| 1 | Von Mises | [4] | B | 245 (+191%) | 3 | 66-81 | 20 | 8.1 |
| 2 | Savage | [5] | B | 309 (+267%) | 4 | 60-62 | 3 | .97 |
| 3 | Reichenbach | [6] | B | 489 (+482%) | 9 | 366-386 | 21 | 4.2 |
| 4 | Venn | [7] | B | 508 (+504%) | 6 ; 10 | 119-166; 235-264 | 78 | 15.3 |
| 5 | Keynes | [8] | B | 539 (+541%) | 7; 8 | 102-124 | 23 | 4.2 |
| 6 | De Finetti | [9] | B | 769 (+815%) | 1; 6; 12 | 10-15; 270-284; 614-626; | 34 | 4.4 |
| 7 | Ramsey | [10] | E | 43 (N.A.) | 1; 2 | 158-166 | 8 | 18.6 |

Table 1 – Textual features of the surveyed works
[Location indicates the first and last pages of each chapter]

Experts of probability are not so familiar with textual analysis and this is the first attempt to break down texts of probability into their components at the best of our knowledge. Textual analysis is not a bibliographical analysis. The latter usually consists in examining the largest number of works dealing with a certain topic; it can be catalogued as an intellectual and subjective survey. The former focuses on a predefined set of works, and requires the researcher to closely investigate the objective content of each work [11]; it can be seen as a statistical and objective inquiry. For example, a researcher of textual analysis can count the number of times certain phrases or words are used in the text; he can define the structure of the text on the basis of the chapters and the sections that compose it; he can dissect author's narrative technique etc.
Let us see the principal outcomes from the textual analysis conducted on the works mentioned above.

 i )    Some sections of the sampled texts are labeled as follows: 'A problem of terminology' [4 p.93], 'The nihilists' [4 p.97], 'The world and the state of the world' [5 p.8], 'The value of observation' [5 p. 125], 'Distinction between logical and psychological view' [7 p.129], 'The application of probability to conduct' [8 p.351], and 'Tyranny of language' [9 p.28]. The titles reveal the concern of writers about qualitative and philosophical themes. The writers argue over on a variety of topics which are distant from mathematics.
 ii )    The authors fill the pages with comments, reflections, notes, explanations, and critical remarks; they share the verbose style of humanists. We use the book of Kolmogorov [3] as comparison term to qualify this aspect of the sample. In fact, Kolmogorov's book does not devote any space to intellectual ruminations and does not make remarks about other views of probability. The percentage increase of pages related to the Kolmogorov book varies from +191% to +815%. The percentage increase cannot be calculated for Ramsey who prepared an essay and not a book.
 iii )    Each master disseminates critical annotations against the concurrent studies. He means to assess his view as the authentic and unique model of probability, but does not demonstrates his perspective using a mathematical proof. In particular each author devotes one or more chapters to disapproving remarks. The right side of Table 1 exhibits the details of this criticism, while Appendix A makes a very short summary of the chapter contents. It may be said that Appendix A expands each line of Table 1 (right side). The rightmost column gives the percentage extent of the chapters respect to the overall work extent and should give an idea of the authors' polemist efforts.
 iv )    Basically, the works have this logical structure: the writer argues about the multifold nature of probability and other arguments too (See point i). He criticizes some models and picks up the model of probability that he judges to be the best in accordance to his proper criteria, finally the author sets up the mathematical construct about the preferred interpretation of *P*. When two authors share a common view, they can even reach opposite conclusions because of their personal approach. For example, Reichenbach [6 p.13] and De Finetti [12] agree that the probability of a single event has no relation with the physical world. The former concludes that this form of probability is to be rejected because of its unreality, the latter places it at the core of his construction.

In conclusion, the textual analysis brings evidence on how the sampled authors – no matter they are frequentist or subjectivist – share a similar style which is strongly based on intellectual discussion, while mathematics ranks second. This method of study did not clarify the foundational issues of probability and statistics so far.

## 3. FOUR OPEN PROBLEMS
Let us recall four significant questions that are still open.

*a) The Problem of Interpretation* – The probability calculus is capable of solving intricate problems and the statistical methods offer support to address sophisticated previsions. Sometime an expert obtains the value of probability *P* at the end of admirable efforts but he is unable to explain whether that number qualifies a material fact or expresses personal credence in that fact. He cannot ensure whether the number *P* is a chance, a possibility, or a wish and the *interpretation problem* is still a debated argument as we have recalled in the introduction of this paper. This fault turns out to be not negligible since statistics has infiltrated several sectors of the present global society and people exploit probability calculations in various areas of modern economies.

*b) The Problem of Discrepancy* – The frequentist and subjective models underpin the classical and the Bayesian statistics in a way [13] and therefore statistical applications should falsify one of the two theories. If the probability models were irreconcilable, as some masters hold, the use of statistics should lend support to one model and should prove the other is false. Instead, countless professional cases demonstrate that both classical and Bayesian statistics are correct and sometimes they furnish identical results. We highlight that a blatant contradiction emerges between theory and practice in the probability domain, and this *discrepancy problem* is awaiting an answer.

*c) The Problem of Choice* – Classical and Bayesian statistics are normally used in professional practice, and sometimes working statisticians wonder: What is the most appropriate statistics to employ in a scheduled project? There is a certain interplay between the two statistics and a precise criterion for selecting the better is missing. An investor who pays for a statistical study and requests for the optimal procedure to follow, does not obtain a precise answer. The rigorous rule to decide on the most appropriate way forward for a planned project is lacking [14], [15] and [16]. Two statisticians may well disagree about the most suitable statistics for given prerequisites; and the *choice problem* is still debated.

*d) The Problem of Axiomatization* – In the year 1900, David Hilbert illustrated a list of problems, unsolved at that time, at the International Congress of Mathematicians held in Paris. The sixth statement of Hilbert suggested the formal exposition of the axioms of the probability calculus [17]. During the 20th century the studies on probability progressed and now we have different axiomatizations in literature. There are also several 'pluralist' writers who accept more than one interpretation [18]; we mention Jean–Antoine Condorcet, Joseph Louis Bertrand, Henry Poincaré, Antoine Cournot, Denis Poisson, Bernard Bolzano, Robert Leslie Ellis, Jacob Friedrich Fries and Karl Popper. Faced with such a plethora of eminent pluralist authors, one might wonder whether Hilbert's sixth problem has several solutions instead of only one; hence one wonders: Is the *axiomatization problem* an ill posed question? The doubts about Hilbert's sixth problem are still in force.

Presently science progresses at high speed; but after three centuries, significant theoretical and practical issues continue to lie at the basis of the probability calculus and statistics.

## 4. THE MATHEMATICAL APPROACH TO DEFINE A MEASURE
Mathematics ranks second in the method adopted by the surveyed authors, instead we mean to employ a method where mathematics ranks first. This is not new. Hilbert in his famous paper wrote: "The investigations on the foundations of geometry suggest the problem: To treat in the same manner, by means of axioms, those physical sciences in which already today mathematics plays an important part; in the first rank are the theory of probabilities and mechanics". Hilbert encouraged theorists to establish the fundamental aspects of probability without any other tool than mathematics. Kolmogorov applies this approach [19] and does not start with philosophical considerations neither makes remarks about other views of probability. Kolmogorov's book presents six axioms and treats purely mathematical topics such as infinite probability fields, random variables, mathematical expectations, conditional probability and independence. We share this approach guided by the mathematical logic since the efforts analyzed in Section 2 left the *Problem of Interpretation* unresolved and in addition have raised or aggravated the *Problems of Discrepancy*, of *Choice* and *Axiomatization*. We decided to address the questions from *a* to *d* using the approach which is recognized as essentially mathematical.

It is necessary to specify that the attribute 'mathematical' pertains to the mode employed by all the physicists, engineers and researchers who determine the new parameter or new measure *y* using mathematical tools while verbal comments and remarks have a secondary position. Obviously explanations are to be added to formal expressions but equations and theorems rank first. This method is a typical part of science theorization and centers on the analytical equation that determines *y* in any aspect – e.g. $y = f(x)$ – more precisely the *argument* or *free variable x* yields the numerical value and also the significance of *y*. The measure *y* can have the *general-abstract* meaning and various *specific-practical* meanings, and this grouping corresponds to the usual subdivision of the mathematical calculus which includes two areas with different properties and scopes [20]. On one side, pure mathematicians treat topics that have a generic relationship with physical reality; on the other side technicians, physicists, economists and other professionals treat topics that have specific relationships with practical applications.

Our personal efforts to apply the formal method produced the book [21] which discusses the importance and the roles that the probability argument has from the theoretical viewpoint. The article [22] published in *Physics Essays* summarizes these initial results. Later the book [23] has calculated the frequentist and subjective probabilities using two theorems. This second part of the inquiry will be summarized in the next pages.

**5. ADVANTAGES OF THE MATHEMATICAL APPROACH IN THE PROBABILITY DOMAIN**
The mathematical approach needs the variable *x* to fix *y*, so it is necessary to establish the argument in advance of defining the probability. Pascal inaugurated the calculus of chance in around 1654 but the precise argument of *P* remained undefined until 1933 when Kolmogorov first fixed it in formal terms [24]. Presently broad literature recognizes that *the random event is the argument of probability* while probability is *the measure of how likely an event will happen*

$$P = P(random\ event). \tag{1}$$

Kolmogorov defines the random event apart from concrete existence and establishes that *E* is a *subset* of the event space Φ

$$P = P(E), \qquad E \in \Phi;\ P \in \mathrm{P}. \tag{2}$$

Thus *P(E)* has *general significance;* it is *the probability in abstract* typical of the axiomatic theory.

In applications, there are two noteworthy arguments that are *the long-term event* $E_n$ – also called *collective* by von Mises and *series* by Venn – and *the single event* $E_1$. We look into the properties of $P(E_n)$ and $P(E_1)$ using a pair of theorems, in particular, the theorems – demonstrated in [23] – illustrate the relationships that exist between the relative frequency and the probabilities $P(E_n)$ and $P(E_1)$ in the order.

*Result #1*: It is assumed the Bernoulli scheme and the concept of probability refers systematically to the probability of success. The former is the well-known *Theorem* or *Law of Large Numbers* (TLN), which in the strong form can be expressed in this way

$$F(E_n) \xrightarrow{a.s.} P(E_n), \qquad as\ n \to \infty. \tag{3}$$

The *Theorem of a Single Number* (TSN) holds that the relative frequency of success is not equal to the probability in a single trial

$$F(E_1) \ne P(E_1), \qquad n = 1. \tag{4}$$

In words, *the frequency gets close to the probability of the long-term event* $P(E_n)$ when *n* tends to infinity; and *the frequency does not match with the probability of the single event* $P(E_1)$. Appendix B includes a lemma that completes the second result and establishes the upper bound for the property (4). For the sake of simplicity, this paper focuses on the single event instead of the set of events defined by the lemma.

Speaking in general, it is the job of science to collect plausible explanations and to use scientific experiments to filter through them, retaining ideas that are supported by the evidence and discarding the others. Testing is at the core of the scientific method. If experiments do not fit with *y* or *y* is not testable in absolute, then researchers conclude that the parameter *y* is not extant in the physical world; it has no actual substance.

TLN and TSN explain how *P* can be validated in the physical reality. Specifically, TLN demonstrates that – at least in principle – the probability of repeated events can be tested, hence $P(E_n)$ is a parameter that exists in the world, it is *an authentic physical quantity*. In addition, assumption $n \to \infty$ is consistent with the classical statistical inference that makes propositions about a population, using data drawn from the population with some form of sampling. This is the first part of the answer to *Problem a,* the second part is more complex.

TSN proves that one cannot test the probability of a single random event. It is not a question of tools or environmental constraints; the theorem proves that never and ever one can control $P(E_1)$ and hence $P(E_1)$ *does not qualify a physical quantity*. This conclusion drawn from TSN fits perfectly with the famous aphorism of de Finetti, "Probability does not exist," which has raised much discussion. Some commentators object that this aphorism sounds like a 'radical' statement, others judge it inappropriate for applications, and so on. In fact, when one expresses the personal opinion *X*, others have a legitimate right to contradict it. Instead, when a theorem proves *X*, either one disproves the theorem or private judgments are not allowed about *X*. In consequence of TSN, *the probability of a single trial must be discarded from the scientific realm and scholars must refuse to accept it.*

*Result #2*: Now a contradiction seems to emerge between the theorem of a single number and the professional needs: scientists should reject $P(E_1)$ but instead they are called for the calculation of $P(E_1)$ every day. How justify this apparent discrepancy between theory and practice?

The answer may be found in the semiotic studies. *Semiotics* teaches us that words and numbers are *signs* or *pieces of information* [25]. By definition a piece of information stands for something. For example, the number 0.32 [= $P(E_n)$] represents a physical quantity. The number 0.32 [= $P(E_1)$] does not represent a physical quantity but is still an item of information. The probability of a single trial has semantic value and is a meaning conveyor. As a consequence of the informational position held by $P(E_1)$, theorists and professionals are allowed to recycle $P(E_1)$ that they should refuse to use. Subjective theorists ascribe a subjective meaning to $P(E_1)$, that is employed to qualify *a credence about the occurrence of* $E_1$. Note how this justification is grounded on semiotic concepts which are amply shared in literature. The present method of study legitimates the subjective probability with the support of the undeniable semantic properties of $P(E_1)$.

The Bayesian inference uses *priors* $P(\theta)$ which can be determined from past information, such as preceding experiments and even using other techniques. A prior is *informative* when it expresses specific, definite information about a variable; it is *uninformative* or *diffuse* if it expresses vague information about a variable such as "all the outcome are equally likely". The posterior probability $P(\theta/x)$ is the probability of the parameters $\theta$ after the observations of *x* with likelihood $P(x/\theta)$

$$P(\theta/x) = \frac{P(x/\theta)P(\theta)}{P(x)}.$$

The informational nature of the subjective probability is consistent with the Bayesian logic which employs prior and actual information through appropriate and non-arbitrary criteria.

In conclusion, *Results 1* and *2* provide two distinct answers to *the interpretation problem* which derive from TLN, TSN and semiotics, and not from personal opinions.

*Result #3*: The probabilities $P(E_n)$ and $P(E_1)$ qualify the very different events $E_n$ and $E_1$ that (3) and (4) specify this way

$$\begin{aligned} n &\to \infty, \\ n &= 1. \end{aligned} \tag{5}$$

These constraints are disjoined in point of logic

$$(n > 1) \; OR \; (n = 1). \tag{6}$$

Therefore, the probability of the long-term event and the probability of a single event do not overlap

$$(n > 1) \; OR \; (n = 1) \;\Rightarrow\; P(E_n) \; OR \; P(E_1). \tag{7}$$

This statement proves that $P(E_n)$ and $P(E_1)$ occur in different circumstances and are not irreconcilable. The two forms of probability apply to situations that do not interfere.

Several researchers claim that the frequentist (or the subjective) model is universal, as long as those researchers ignore assumptions (5) that are typical of the theorems of large numbers and a single number. The present frame shows how $P(E_n)$ and $P(E_1)$ regard specific practical situations, while $P(E)$ has generic relationship with the physical reality. Each form of applied probability is subjected to special experimental constraints and should not be confused with the abstract form (2). The confusion between applied and abstract forms of probability causes the *discrepancy problem* whereas mathematical statement (7) proves this problem does not have a basis.

*Result #4*: The theorems of large numbers and a single number take mutually exclusive hypotheses (5) thus the statistical methods underpinned by TLN and TSN must be consistent with (5). The logical divide (6) yields a second divide that is expressed as follow:

> *If one means to investigate a long-term event, he must resort to using classical statistics;*
> *If one means to focus on a single event, he must adopt Bayesian statistics.*

Statisticians employ an assortment of mathematical tools; some of these are frequentist, and some are Bayesian. There are Bayesian methods which have a frequentist equivalent and sometimes give the same numerical result. In practice, there are situations in which one of the methods is more preferred by some criteria, while the other method is preferred for other reasons. There are also scholars who believe that each statistics is actually essential for full development of the other. Bayarri and Berger tell that a certain interplay occurs between classical and Bayesian statistics [26]. Someone even holds a hybrid approach [27]. By contrast, the present frame shows two distinct ways: assumption $n \to \infty$ is consistent with the classical statistics, and $n = 1$ is consistent with the Bayesianism. Statements (8) say something more and above; they provide a rule that does not admit exceptions. There is no middle way in (8) and this rule does not allow exceptions. For example, if a doctor is treating a patient who is suffering from cancer Z, he is considering a single case and adopts the Bayesian tools. If a researcher investigates the epidemiology of cancer Z, he is concerned with a general trend and follows classical statistics. The working statistician has to apply a precise statistical technique depending on the particular question he is dealing with. The two statistical schools work perfectly in practice provided that they are restricted to a suitable domain of application and the rule (8) constitutes the answer to the *Problem of Choice*.

The subjective model underpins the Bayesian statistics, which focuses on single occurrences even if the Bayesian procedures are not confined to a lone observation. When a Bayesian applies to a sequel of repeated events, his conclusions regard each individual case. Results #1 and #2 explain how the significance of the frequency probability does not have anything in common with the subjective meaning, and thus each personal value $P(E_{11})$, $P(E_{12})$, $P(E_{13})$,… $P(E_{1n})$ cannot be confused with $P(E_n)$ that is an authentic physical parameter, regardless of the numerical values which may be identical.

*Result #5*: The mathematical approach which we are employing leads to the precise organization that divides the abstract calculus from the applied calculus. The probability sector splits into two areas in consequence of the arguments $E$, and $E_n$ with $E_1$. The generalized definition of $P(E)$ is unique in accordance with the expectations of Hilbert's sixth problem that proves not to be an ill posed question (*Problem of axiomatization*). The applied calculus regards special cases such as the long-term event $E_n$, the single event $E_1$, and even the quantum event $E_q$, the economic event $E_x$ etc. There are various areas of application, while the abstract theory is single.

This precise organization of the probability calculus does not falsify the mathematical constructions of the frequentist and subjective authors, rather it shows how those constructs are incomplete. Each theory provides the illustration of the probability calculus under the explicit hypothesis $n \to \infty$ and the hypothesis $n = 1$ in the order, hence they offer effective assistance to scientists in each application field, but do not provide the exhaustive illustration of the probability concept.

The reader can note how the method that we have adopted provide answers to the *Problems* from *a* to *d* which are marked by innovation.

## 7. CONCLUSION AND OUTLOOK

Three hundred and fifty years after Pascal's inauguration of the modern probability calculus there are significant open problems. The textual analysis provides evidence that opinions and personal choices often influence the works of masters. They not only leave unresolved some foundational issues but sometimes compound them. At the present time, experts struggle with the *discrepancy problem*, the *choice problem* and the *axiomatization problem* besides the *interpretation problem*. This failure pushed us to follow the way suggested by Hilbert and inaugurated by Kolmogorov which centers on mathematics while verbal explanations serve as completing elements. This approach provides innovative answers to the *Problems a-d*; more precisely:
- It is proved the *testability* of $P(E_n)$ and the *unverifiability* of $P(E_1)$;
- It is explained how *$P(E_1)$ is reused as subjective probability* instead of being scrapped;
- It is demonstrated why *$P(E_n)$ and $P(E_1)$ are not irreconcilable*;
- *Statisticians have a precise rule to follow* when they are called for selecting the most appropriate statistics;
- The probability calculus splits into the abstract calculus of $P(E)$ and the applied calculus of $P(E_n)$ and $P(E_1)$, hence *the Hilbert's sixth problem is not ill posed*.

Several experts of probability and statistics are inclined to accept both the frequentist and subjective models (see Section 3) but a rigorous theoretical frame is missing so far. To the best of our knowledge, this inquiry provides the first formal 'dualist theory'.

Blanshard [28] claims that *philosophy comes before science,* especially when a problem is not well-defined. It is difficult to conceive of a measure when one does not quite know what he is measuring, however the importance and role of philosophy change with the progress of science. When scientists find the mathematical definition of the searched measure and this equation reaches broad consensus then the intellectual ruminations progressively slip to the background. For instance, thinkers debated the nature of mechanical force for a long while. When Newton fixed the concept of force with a differential equation, engineers began to calculate mechanical equipment and verbose discussion became unfashionable. The debates about the multifold nature of probability were appropriate in the past since the concept of probability could not be accurately stated, the context of the problem seemed fuzzy. Nowadays the theorems of large numbers and a single number demonstrate the properties of probability with precision and physical phenomena wait to be calculated using the present results as first in quantum mechanics. In particular we are going to present a study that casts new light into the quantum wave collapse.

**APPENDIX A – A SHORT REVIEW**

**1.** Richard von Mises devotes the initial part of the third chapter in the book [4] to investigating the opposing definitions of probability; the second part of the third chapter deals with the critical aspects of his theory. Von Mises emphasizes the limit of the Laplacian scheme, and specifically holds that the subjective model may be influenced by psychological or physiological mechanisms and puts down this model as follows: "The peculiar approach of subjectivists lies in the fact that they consider 'I presume that these cases are equally probable' to be equivalent to 'These cases are equally probable' since for them probability is only a subjective notion".

**2.** Savage writes three pages in Chapter 4 of *The Foundations of Statistics* [5] with critical remarks on the 'objective' interpretations of probability. His view may be summed up with this passage: "In the first place, objectivistic views typically attach probability only to very special events. Thus on no ordinary objectivistic view would it be meaningful... Secondly objectivistic views are .... charged with circularity. They are generally predicated on the existence in nature of processes that may ... be represented by ... an infinite sequence of independent events."

**3.** Reichenbach illustrates the historical evolution of the probability concept in the first chapter of [6]. Chapter nine discusses the various meanings of probability; in particular, he examines the probability of a single event, which should have no place in science. If the concept of probability only represents a subjective expectation, then *P* does not have any connection to the real world.

**4.** John Venn criticizes the interpretation of probability as personal belief – especially in relation to the thought of De Morgan – in Chapter VI of the book [7]. He states: "the difficulty of obtaining any measure of the amount of our belief" and adds: "we experience hope or fear in so many very instances, that … whilst we profess to consider the whole quantity of our belief we will in reality consider only a portion of it." Venn concludes that human actual belief is one of the most elusive and variable factors so that we can scarcely ever get sufficiently clear hold of it to measure it. Chapter X tackles another argument, specifically it questions whether the events calculated with the probability calculus are to be attributed to chance on one hand or alternatively to causation or design on the other hand.

**5.** Chapter VII of Keynes's treatise [8] makes an historical retrospect of the probability calculus. He illustrates frequency theory recalling the work of Leslie Ellis and mostly looking into the Venn work. Chapter VIII emphasizes the limits of the frequentist model that clearly excludes a great number of judgments which are generally believed to deal with probability. Keynes also stresses the practical use of statistical frequencies, since "an event may possess more than one frequency, and that we must decide which of these to prefer on extraneous grounds." Later Keynes emphasizes the differences between Venn's construction and the generalised frequency theory which he means to put forward.

**6.** De Finetti illustrates his theory in two volumes [9] that are peppered with unfavorable and even sarcastic judgments about the opponent theories. In the first chapter he places the notions typical of the subjective and objective schools of probability side by side in order to highlight the profound differences extant between them. For instance, he notes that for 'objectivists' two events are independent if the occurrence of one does not affect the probability of the other; instead for 'subjectivists' two events are independent if the knowledge of one does not modify the assessment of the probability of the other event. Chapter six introduces three main interpretations of distributions, then the author begins a long discussion against countable additivity. Chapter twelve deals with estimations and testing that have distinct characters from the perspectives of the classical and Bayesian statistics. De Finetti never fails to emphasize the Bayesian techniques and to criticize the alternative methods.

**7.** Ramsey begins the seventh chapter of the essay [10] with censorious comments on the works of von Mises and Keynes. He pinpoints that the latter recognizes the subjectivity of probability but in substance Keynes does not assign any value to subjectivism. Moreover, Keynes believes there is an objective relationship between knowledge and probability, as knowledge is disembodied and not personal. Ramsey analyses the connection between the subjective degree of belief an individual has in a proposition and the probability it can be given. As regards the frequentist theory he writes: "I am willing for the present to concede to the frequency theory that probability as used in modern science is really the same as frequency."

**APPENDIX B – UPPER-BOUND LEMMA**
TSN proves that probability is unreal when the argument is a single event. We could say that *n*=1 is the lower bound of the probability inexistence. Let us examine the largest number of events whose probability is unreal.

*Lemma*: Suppose *z* is any positive integer and the probability of *E* verifies

$$P(E) = 1/z, \qquad z > 0. \tag{B.1}$$

Then the relative frequency of the successful event *E* in *n* trials is not equal to the probability

$$F(E_n) \neq P(E). \tag{B.2}$$

If

$$1 < n < z. \tag{B.3}$$

*Proof:* We proceed by absurd and deny (B.2), we put the relative frequency $F(E_n) = N(E_n)/n$ equals to the probability

$$N(E_n)/n = P(E). \tag{B.4}$$

When the event *E* occurs one time

$$N(E_n) = 1.$$

We obtain from (B.4) and (B.1)

$$1/n = 1/z.$$

And thus

$$n = z.$$

This conclusion mismatches with assumption (B.3), therefore (B.4) is false and (B.2) true.
*Example*: The probability of getting a king from a card deck

$$P(E_K) = 4/52 = 1/13 = 1/z.$$

If the number of trials is less than thirteen

$$13 > n > 1. \tag{B.5}$$

There are two possibilities. If one does not get any king in *n* drawings, the relative frequency is lower than $P(E_K)$

$$0/n < 1/13.$$

If one gets one (or more) kings in *n* experiments, the relative frequency is greater than $P(E_K)$. Suppose to minimize the number of success and minimize the number of drawings we get

$$1/12 > 1/13$$

In summary, the relative frequency never collides with the probability $P(E_K)$ that proves to be out-of-control and therefore can but express a personal degree of belief in relation to the condition (B.5).

**References**

[1] I. Hacking, *The Emergence of Probability: A Philosophical Study of Early Ideas About Probability Induction and Statistical Inference* (Cambridge University Press, 1975).
[2] T. Fine, *Theories of Probability* (Academic Press, 1973).
[3] A. Kolmogorov, *Foundations of the Theory of Probability* (Chelsea, 1950).
[4] R. von Mises, *Probability, Statistics and Truth* (Mc Millan Co., 1957).
[5] L.J. Savage, *The Foundations of Statistics* (Wiley, 1954).
[6] H. Reichenbach, *The Theory of Probability: An Inquiry into the Logical and Mathematical Foundations of the Calculus of Probability* (University of California Press, 1949).
[7] J. Venn, *The Logic of Chance* (McMillan Co., 1888).
[8] J.M. Keynes, *A treatise of Probability* (McMillan & Co., 1921).
[9] B. de Finetti, *Teoria della Probabilità* (Einaudi, 1970). Translated as *Theory of Probability* (John Wiley and Sons, 1974).
[10] F.P. Ramsey, Truth and Probability, In *The Foundations of Mathematics and other Logical Essays*, Braithwaite R.B. (ed), (Routledge & Kegan Paul, 1931): 156-198.
[11] B. Helder, *Textual Analysis: An Approach to Analysing Professional Texts* (Samfundslitterature, 2011).
[12] M.C. Galavotti, Anti-realism in the Philosophy of Probability: Bruno de Finetti's Subjectivism, *Erkenntnis*, 31: 239-261(1989).
[13] D.R. Cox, Statistics, an Overview, In *Encyclopedia of Biostatistics* (John Wiley & Sons, 2005) 7.
[14] S.C Chow., *Controversial Statistical Issues in Clinical Trials* (Chapman & Hall/CRC, 2011)
[15] R.E. Weiss, Bayesian Methods for Data Analysis, *American J. of Ophthalmology,* 149(2): 187-188 (2010).
[16] K.J. Friston, W. Penny, C. Phillips, S. Kiebel, G. Hinton, J. Ashburner, Classical and Bayesian Inference in Neuroimaging: Theory, *NeuroImage,*16: 465-483 (2002).
[17] R. Thiele, Hilbert and his twenty-four problems, in *Mathematics and the Historian's Craft* (Springer, 2005): 243-295.
[18] W.C. Salmon, Dynamic rationality, In Fetzer J. (ed), *Probability and Causality* (Reidel Publishing Company,1988).
[19] G. Shafer, V. Vovk, The Sources of Kolmogorov's Grundbegriffe, *Statistical Science,* 21(1): 70–98 (2006).
[20] C. Pincock, Towards a Philosophy of Applied Mathematics, In *New Waves in Philosophy of Mathematics*, Bueno O. and Linnebo Ø. (eds) (Palgrave Macmillan, 2009): 173-194.
[21] P. Rocchi, *Structural Theory of Probability*, Springer (2003).
[22] P. Rocchi, L. Gianfagna, Probabilistic Events and Physical Reality: A Complete Algebra of Probability, *Physics Essays*, 15(3): 331-118 (2002).
[23] P. Rocchi, *Janus-Faced Probability* (Springer, 2014).
[24] P. Rocchi, De Pascal à Nos Jours: Quelques Notes sur l'argument *A* de la probabilité *P*(*A*), *Actes du Congrès Annuel de la Société Canadienne d'Histoire et de Philosophie des Mathématiques,* vol. 19: 228-235 (2006).
[25] D. Chandler, *Semiotics: The Basics* (Routledge, third edition, 2017).
[26] M.J. Bayarri, J.O. Berger, The Interplay of Bayesian and Frequentist Analysis, *Statistical Science*, 19(1):58-80 (2004).
[27] A. Yuan, Bayesian Frequentist Hybrid Inference, *Ann. Statist.*, 37(5A): 2458-2501 (2009).
[28] B. Blanshard, The Philosophic Enterprise, In *The Owl of Minerva: Philosophers on Philosophy*, C.J. Bontempo and S.J. Odell (eds) (McGraw-Hill, 1975): 163-177.